%% file: REDUCTIBLE.tex
\documentclass[a4paper,11pt]{amsart}
\usepackage{latexsym, delarray}
\usepackage[all]{xy}
\usepackage{amssymb}

\usepackage[colorlinks=true]{hyperref}

\newcommand{\gal}{\mathrm{Gal}}
\newcommand{\qbar}{\overline{\mathbb{Q}}}
\newcommand{\Frob}{\mathrm{Frob}}
\newcommand{\SL}{\mathrm{SL}}
\newcommand{\GL}{\mathrm{GL}}

\newcommand{\Flbar}[1][l]{\overline{\mathbb{F}}_{#1}}
\newcommand{\CuspF}[3]{\mathcal{S}_{#1}(\Gamma_{1}(#2),#3)}
\newcommand{\CuspFsans}[2]{\mathcal{S}_{#1}(\Gamma_{1}(#2))}
\newcommand{\CuspFtriv}[2]{\mathcal{S}_{#1}(\Gamma_{0}(#2))}
\newcommand{\tr}{\mathrm{tr}}
\newcommand{\zbar}{\overline{\mathbb{Z}}}

\newcommand{\eps}{\varepsilon}

\newcommand{\FF}{\mathbb{F}}

\newcommand{\ce}{\mathbb{C}}
\newcommand{\erre}{\mathbb{R}}
\newcommand{\z}{\mathbb{Z}}
\newcommand{\ene}{\mathbb{N}}
\newcommand{\q}{\mathbb{Q}}

\newcommand{\G}{\Gamma_0}

\newcommand{\trivial}{\mathbf{1}}

\newcommand{\rep}{\trivial\oplus\chi_l^{k-1}}

\newcommand{\Id}{\mathrm{Id}}
\newcommand{\lcm}{\mathrm{lcm}}

\newtheorem{thme}{Theorem}
\newtheorem{thm}{Theorem}[section]
\newtheorem{lema}[thm]{Lemma}
\newtheorem{prop}[thm]{Proposition}
\newtheorem{cor}[thm]{Corollary}
\newtheorem{conj}[thm]{Conjecture}

\theoremstyle{definition}
\newtheorem{rem}[thm]{Remark}

\begin{document}

\title[Modularity of reducible Galois representations]{On the modularity of reducible \\ mod~$l$ Galois representations}
\author[Nicolas Billerey]{Nicolas Billerey $^\P$}
\thanks{\P Universit\'e Blaise Pascal Clermont-Ferrand~2}
\address{Laboratoire de Math\'ematiques, Universit\'e Blaise Pascal Clermont-Ferrand~2, Campus universitaire des C\'ezeaux, 63177 Aubi\`ere Cedex, France}
\email{Nicolas.Billerey@math.univ-bpclermont.fr}
\author[Ricardo Menares]{Ricardo Menares $^\dag$ }
\address{Instituto de Matem\'aticas, Pontificia Universidad Cat\'olica de Valpara\'iso, Blanco Viel 596, Cerro Bar\'on, Valpara\'iso, Chile}
\thanks{\dag Pontificia Universidad Cat\'olica de Valpara\'iso. Partially supported by FONDECYT grant 11110225 and CONICYT grant Inserci\'on en la academia}
\email{ricardo.menares@ucv.cl}


\maketitle

\begin{abstract}
Given an odd, semisimple, reducible, $2$-dimen\-sional mod~$l$ Galois representation, we investigate the possible levels of the modular forms giving rise to it. When the representation is the direct sum of the trivial character and a power of the mod~$l$ cyclotomic character, we are able to characterize the  primes that can arise as levels of the associated newforms. As an application, we determine a new  explicit lower bound for the highest degree among the fields of coefficients of newforms of trivial Nebentypus and prime level. The bound is valid in a subset of the  primes with natural (lower) density at least 3/4. 
\end{abstract}

\input{intro2}
\input{eisenstein}

\input{modularitereductible}

\input{modularitecyclotomique}

\input{corps}


\end{document}

%% file: intro2.tex
\section*{Introduction}

Let $l$ be a rational prime number. In this paper we are interested in $2$-dimensional mod~$l$ Galois representations, that is, continuous homomorphisms

\begin{equation}\label{representacion}
 \rho~:\gal(\qbar/\q)\longrightarrow\GL_2(\Flbar).
\end{equation}

Work of Deligne, extending earlier results by Eichler and Shimura, shows that such representations arise naturally in the theory of modular forms. More precisely, let us take $\qbar$ to be the algebraic closure of~$\q$ in~$\ce$ and fix a place~$w$ of~$\qbar$ over~$l$. We denote by $a\mapsto\widetilde{a}$ the reduction map modulo $w$ from the ring of integers~$\zbar$ of $\qbar$ to the residue field~$\Flbar$. Let us denote by~$\CuspFsans{k}{N}$ the $\ce$-vector space of cuspidal modular forms of weight~$k\ge2$ for~$\Gamma_1(N)$. Then, attached to any form $f\in\CuspFsans{k}{N}$ that is an eigenform for the Hecke operators~$\{T_p\}_{p\nmid N}$ with corresponding set of eigenvalues~$\{a_p\}_{p\nmid N}$, there is a Dirichlet character $\chi$ of modulus~$N$ and an odd semisimple Galois representation   
\[
\rho_f~:\gal(\qbar/\q)\longrightarrow\GL_2(\Flbar),
\]
unique up to isomorphism, that is unramified outside~$Nl$ and satisfies for every prime~$p\nmid Nl$
\[
\left\{
\begin{array}{l}
\tr(\rho_f(\Frob_p)) = \widetilde{a_p} \\
\det(\rho_f(\Frob_p)) = \widetilde{\chi(p)}p^{k-1} \\
\end{array}
\right.
\]
where $\Frob_p$ denotes a Frobenius element at~$p$ in~$\gal(\qbar/\q)$.

A natural problem is then to understand which mod~$l$ Galois representations~$\rho$ are modular, meaning that $\rho$ is isomorphic to some~$\rho_f$ as above. We shall say in this case that~$\rho$ arises from the cuspidal eigenform~$f$.  Since $\rho$ does not uniquely determine $f$, it is  equally important to determine the set of  possible values of  the level, weight and character of $f$. In this article, we will focus on the level aspect. 

In the case where $\rho$ is \emph{irreducible}, these questions are addressed  by Serre's modularity conjecture  (\cite{Ser87}), which is  now a theorem of Khare and Wintenberger (\cite{KhWi09a,KhWi09b}). Thanks to these works (building on seve\-ral deep results by many people), it is known that  every odd, irreducible Galois representation~$\rho$ as in~\eqref{representacion} arises from a cuspidal eigenform. Let $N(\rho)$ be the prime-to-$l$ part of the Artin conductor of~$\rho$. Carayol (\cite{Car89}) and Livn\'e (\cite{Liv89}) independently proved that $N(\rho)$ is minimal (for divisibility) among all possible prime-to-$l$ levels of modular forms giving rise to~$\rho$. Ribet's famous level-lowering theorem (\cite[Thm.~1.1]{Rib90}) shows that a modular~$\rho$ indeed arises from a Hecke eigenform of this `optimal' level. On the other hand, Diamond and Taylor soon thereafter (\cite{DiTa94}) completely described the set of prime-to-$l$ integers~$M>N(\rho)$ such that $\rho$ arises from a weight-$k$ {\em newform} of level~$M$, in terms of the ramification theory of $\rho$. They called such integers  \emph{non-optimal levels} attached to~$\rho$.


In this paper we focus on the case where $\rho$ is \emph{reducible}. Given an even integer~\(k\geq2\), we have that, if $l\ge k-1$, then every modular mod~$l$ reducible representation arising from a weight-$k$ newform with squarefree level and trivial Nebentypus is of the form $\trivial\oplus\chi_l^{k-1}$, where $\chi_l$ is the mod~$l$ cyclotomic character (cf. Proposition~\ref{ss:reducible_mod_reps}). The results of this paper are most precise when~$\rho$ is of this form. 

The representation  $\trivial\oplus\chi_l^{k-1}$  is unramified outside~$l$, hence the prime-to-$l$ part of  its Artin conductor is 1. Because of this fact, we will say that $\trivial\oplus\chi_l^{k-1}$ is \emph{modular of optimal level}  if it arises from a cuspidal eigenform of weight $k$ and level 1. Contrary to what happens in the irreducible case, in spite of the fact that $\trivial\oplus\chi_l^{k-1}$ is odd and semi-simple,  this representation   need not be modular of optimal level for every choice of $k$ and $l$. Indeed, when $k=2$, there are no cusp forms of level one and weight two and hence $\trivial\oplus\chi_l$ is not modular of optimal level.  On the other hand,  in (even) weight~$k\geq 4$, Ribet has  proved that, if $l>k+1$, then the representation $\trivial\oplus\chi_l^{k-1}$ is modular of optimal level  if and only if $l$ divides the numerator of~$B_k/k$ where~$B_k$ denotes the $k$-th Bernoulli number (more precisely, the direct implication  of this assertion results from  \cite[Lem.~5.2]{Rib75} and the reverse implication is mentioned in \cite[Prop.~1]{Gha02}).

Our goal is  to study the set of integers $N\ge2$, prime to $l$, such that $\trivial\oplus\chi_l^{k-1}$  arises from a newform of weight $k$, level $N$ and trivial Nebentypus. Following the terminology that Diamond and Taylor introduced  in the irreducible case, we call such integers \emph{non-optimal levels} attached to this representation.

Our main concern is the classification of the squarefree, non-optimal levels attached to $\trivial\oplus\chi_l^{k-1}$.  In weight $k=2$, the first step in this study is a well-known result of Mazur (\cite[Prop.~(5.12)]{Maz77}), asserting that $\trivial\oplus\chi_l$ arises from a weight-$2$ cusp form of prime level~$N$ and trivial Nebentypus if and only if $l$ divides the numerator of~$(N-1)/12$. Then,  Ribet  classified the squarefree non-optimal levels having exactly two prime factors,  under the assumption~$l>3 $ (cf.~\cite{Rib10}). His results have been  very recently extended to the case of three prime factors by his student~Hwajong Yoo in his Ph.D. thesis (Spring 2013). These works show that the classical level-raising condition does not suffice to determine the non-optimal levels.

Our main theorem  treats the case $k \geq 4$ and is an interpolation of Mazur's and Ribet's results quoted above  in the following sense.

\begin{thme}\label{main}
Let $k$ be an even integer~$\ge4$ and assume $l>k+1$. Then the representation~$\trivial\oplus\chi_l^{k-1}$ arises from a weight-$k$ newform of prime level~$N$ and trivial Nebentypus if and only if at least one of the following conditions holds~:
\begin{enumerate}
\item \label{prima} $N^k\equiv 1\pmod{l}$
\item \label{dona} $N^{k-2}\equiv 1\pmod{l}$ and $l$ divides the numerator of $B_k/k$.
\end{enumerate}
\end{thme}

We conjecture that, when $k\geq 4$ is even,  the 
level-raising condition at a prime $p$ is actually sufficient for the representation~$\trivial\oplus\chi_l^{k-1}$  to arise from a newform of squarefree level divisible  by $p$. This leads to  a conjectural description of the set of squarefree levels of weight-$k$ newforms with trivial Nebentypus that give rise to~$\trivial\oplus\chi_l^{k-1}$ (cf. Conjecture~\ref{conj} and the ensuing remarks).  Theorem \ref{main} confirms this description for prime non-optimal levels. We discuss in Section~\ref{ss:level_raising_condition} how hypothesis~\eqref{prima} and~\eqref{dona} relate to the level-raising condition.

We stress the fact that the statement of Theorem~\ref{main} considers not only eigenforms, but newforms (in Mazur's  and Ribet's statements this distinction is unnecessary).  We deduce this statement from explicit computations of the constant term of various Eisenstein series at the cusps of an appropriate modular curve. We carry on  these computations in Section~\ref{s:Eisenstein}, following a classical analytic approach. These calculations refine related constant term computations appearing in the work of  Faltings and Jordan  \cite{FaJo95},  who use algebraic geometry methods through the interpretation of modular forms as sections of a line bundle on a modular curve (cf. Remark \ref{overlap}). Our proof of the reverse implication in Theorem~\ref{main} proceeds by first constructing, using Eisenstein series modulo~$l$ and the Deligne-Serre lemma, a cuspidal eigenform satisfying the necessary congruences. We then use a theorem of Diamond in~\cite{Di91} to show that we can take the eigenform to be a newform of level~$N$.



A related general question is to decide whether a given odd, reducible, but not necessarily semisimple, mod~$l$ Galois representation is the reduction of a $l$-adic representation attached to some cuspidal eigenform. We shall prove in Section~\ref{s:modularite_reductible} that every odd representation which is the sum of two characters arises from a cuspidal modular form (Theorem~\ref{thm:modularity_general}). This fact is presumably well-known to experts, but our methods allow us to  provide an elementary and self-contained proof and we include it here because we could not locate a proof in the literature. We refer the interested reader to the work of Ramakrishna~\cite{Ram02}, and the references therein, for a collection of important results on this question.

We  apply our results to the following situation. To a normalized eigenform  $f \in \CuspFtriv{k}{N}^{\mathrm{new}}$ with $q$-expansion at infinity $$f(z)=\sum_{n \geq 1 }^{ } a_n(f) q^n, \quad q=e^{2\pi i z},$$ we attach the number field  $K_f := \q \big( a_n(f) : (n,N)=1 \big).$ Put 
$$d^{\mathrm{new}}_k(N):= \max \{ [K_f:\q] : f \in \CuspFtriv{k}{N}^{\mathrm{new}}, \textrm{ normalized Hecke eigenform} \}.$$ 

A Theorem of Royer\footnote{Royer's theorem holds for arbitrary levels~$N$ with a constant depending on a fixed prime not dividing~$N$. It is however stated only in the case $k=2$ in \emph{loc. cit.} but the proof undoubtedly extends to weights~$\ge 2$.} (\cite{Ro00}) implies that 
\begin{equation}\label{Royer}
d^{\mathrm{new}}_k(N) \gg_k \sqrt{\log \log N}, \quad N \rightarrow \infty, \quad  N \textrm { prime}.
\end{equation}

As an application of Theorem~\ref{main}, we obtain, in the spirit of~\cite{DiUrRi11}, a lower bound for $d^{\mathrm{new}}_k(N)$. As a new ingredient, we use results from analytic number theory on the density of prime numbers~$p$ for which $p-1$ has a large prime factor (cf. \cite{Gold69}, \cite{Luca}).  We then manage to obtain a bound  which is better than~\eqref{Royer} but is only valid in  a restricted class of prime numbers. 

\begin{thme}\label{thm:degree}
There exists an explicit set  of primes $\mathcal{P}$ of (natural) lower density at least $\frac{3}{4}$ with the property that, for every even integer $k \geq 2$, there exists a constant $c_k>0$   such that the inequality  $$d^{\mathrm{new}}_k(N) \geq c_k\log N $$ holds for all  $N \in \mathcal{P}$ with $N\geq (k+1)^4$. The constant $c_k$ can be taken as  $$c_k=\left(8\log  \left(1+2^{(k-1)/2}\right) \right)^{-1}.$$  
\end{thme}

If we assume the truth of Conjecture \ref{conj}, then it is possible to extend the validity of the above bound to appropriate squarefree integers (cf. Theorem~\ref{souslaconjecture}). 

In the spirit of Maeda's conjecture, Tsaknias has proposed a conjectural lower bound for $d_k^{\textrm{new}}(N)$ for $N$ fixed and varying $k$ \cite{Tsa14}. His conjecture implies that there exists  a constant $c>0$ such that, for all prime numbers $N$, there is  an integer $k(N)$ such that $d_k^{\textrm{new}}(N)>cN$ for all $k \geq k(N)$. Further numerical data, that he has generously shared with us, suggest that $k(N)$ is a bounded function of $N$. If this were true, then~$d_k^{\textrm{new}}(N)$  would grow linearly with $N$ if $k$ is fixed. 

\medskip
\noindent{\bf Acknowledgements. } The authors benefited from the warm hospitality of the mathematics department at  Universidad T\'ecnica Federico Santa Mar\'ia. We also received from  \'E.~Fouvry, F.~Luca and E.~Royer  explanations and concrete suggestions on the analytic number theory surrounding Theorem~\ref{thm:degree}. We thank K. Ribet and G.~Wiese for interesting comments and references.  P.~Tsaknias has kindly  shared numerical  data with us. We are grateful to the anonymous referee, whose careful reading of the manuscript greatly helped  to improve it.

%% file: eisenstein.tex
\section{Preliminaries on Eisenstein series}\label{s:Eisenstein}
In this section we recall some classical definitions and compute the constant term of the $q$-expansion at various cusps of some specific Eisenstein series that will be used in the sequel. 

\subsection{Gauss sums and Bernoulli Numbers}\label{ss:Gauss_sum_Bernoulli_Numbers}
For an integer $m \geq 2$, we set $$C_m=\frac{(-2i\pi)^m}{(m-1)!}.$$ Let $\psi:(\z/c\z)^{\times}\rightarrow\ce^{\times}$ be a primitive Dirichlet character of modulus~$c\ge 1$. The Gauss sum attached to~$\psi$ is defined by
\begin{equation*}
W(\psi)=\sum_{n=1}^{c}\psi(n)e^{2i\pi n/c}
\end{equation*}
and the Bernoulli numbers $(B_{m,\psi})_{m\ge1}$ by
\begin{equation*}
\sum_{n=1}^c\psi(n)\frac{te^{nt}}{e^{ct}-1}=\sum_{m\ge0}B_{m,\psi}\frac{t^m}{m!}.
\end{equation*}
In particular, if $\psi$ is the trivial character (of modulus~$1$), then $B_{m,\psi}$ is the classical Bernoulli number~$B_m$, except when $m=1$ in which case we have $B_{1,\psi}=-B_1=1/2$.

The Bernoulli numbers are related to certain special values of the $L$-function $L(s,\psi)$ attached to~$\psi$. More precisely, we have the following proposition (which follows, for instance, from \cite{Was97}, Theorem~4.2 and  the functional equation on p. 30 of \emph{loc. cit.}).
\begin{prop}\label{prop:FE_for_L_functions}
Let $m\ge 2$ be an integer such that $\psi(-1)=(-1)^m$. Then, we have that
\begin{equation*}
L(m,\psi)=-W(\psi)\frac{C_m}{c^m}\cdot\frac{B_{m,\overline{\psi}}}{2m}\not=0,
\end{equation*}
where $\overline{\psi}$ means the complex conjugate of~$\psi$.
\end{prop}

\subsection{Constant term  computations}
For a positive real number $A$, let us denote by $\alpha_A$ the operator acting on complex valued functions $f$ on the upper half-plane~$\mathfrak{h}$ by
\[
\alpha_A(f)(z)=f(Az).
\]
\subsubsection{General computations}
Let $k\ge3$ be an integer and $\varepsilon_0:(\z/N\z)^{\times}\rightarrow\ce^{\times}$ be a primitive Dirichlet character of modulus~${N\ge1}$ such that $\varepsilon_0(-1)=(-1)^k$. We denote by $E_k^{\trivial,\varepsilon_0}$ the Eisenstein series in~$\mathcal{M}_k(\Gamma_1(N),\varepsilon_0)$ given by the following $q$-expansion~:
\begin{equation}\label{eq:Eisenstein_series_main}
E_k^{\trivial,\varepsilon_0}(z)=-\frac{B_{k,\varepsilon_0}}{2k}+\sum_{n\ge1}\Big(\sum_{m\mid n}\varepsilon_0(m)m^{k-1}\Big)q^n, \quad\text{where } q=e^{2\pi i z}.
\end{equation}
Note that when $N=1$, $E_k^{\trivial,\varepsilon_0}=E_k^{\trivial,\trivial}$ is nothing but the classical level~$1$ Eisenstein series of weight~$k$
\begin{equation}\label{eq:classical_Ek}
E_k=-\frac{B_k}{2k}+\sum_{n\ge 1}\Big(\sum_{m\mid n}m^{k-1}\Big)q^n.
\end{equation}
Note also that our function~\(E_k^{\trivial,\varepsilon_0}\) differs from that of~\cite[Thm.~4.5.1]{DiSh05} by a factor~\(2\) (that is, their function~\(E_k^{\trivial,\varepsilon_0}\) is twice ours).
The main goal of this paragraph is to compute the constant term of the $q$-expansion of $\big(\alpha_ME_k^{\trivial,\varepsilon_0}\big)\vert_k\gamma$ where $M\ge 1$  is an integer coprime to~$N$, $\gamma\in\SL_2(\z)$ and the notation $\vert_k$ means the classical slash operator on modular forms.

\begin{prop}\label{prop:constant_term}
Let $\varepsilon_0$ and $k$ as above. Let $M$ be an integer~$\ge1$ coprime to~$N$ and $\gamma=
\begin{pmatrix}
u & \beta \\
v & \delta \\
\end{pmatrix}\in\mathrm{SL}_2(\z)$. The constant term of the $q$-expansion of~$\big(\alpha_ME_k^{\trivial,\varepsilon_0}\big)\vert_k\gamma$ is
\[
\left\{
\begin{array}{cl}
0 & \textrm{if }\displaystyle{N\nmid\frac{v}{r}} \\
\displaystyle{-\frac{\overline{\varepsilon_0}(M')\varepsilon_0(\delta)}{M'^k}\frac{B_{k,\varepsilon_0}}{2k}\not=0} & \textrm{otherwise}
\end{array}
\right.
\]
where $r=\gcd(v,M)$ and $M'=M/r$.
\end{prop}

\begin{rem}\label{overlap}
The constant term of this kind of Eisenstein series plays a role in the work of Faltings and Jordan (cf.  \cite{FaJo95}, Definition 3.16 and Theorem 3.20). However, they compute these constant terms only up to a unit in an appropriate ring of integers. 
\end{rem}

\proof With the notations of Section \ref{ss:Gauss_sum_Bernoulli_Numbers} put
\[
G=\frac{2C_kW(\overline{\varepsilon_0})}{N^k}E_k^{\trivial,\varepsilon_0}.
\]
According to the remark above, the function~\(G_k^{\trivial,\varepsilon_0}\) of~\cite[Thm.~4.5.1]{DiSh05} is twice our function~\(G\). It then follows from the definition of this function on page~127 of \emph{loc. cit.} that 
\[
G=\sum\limits_{\substack{j=0\\ \gcd(j,N)=1}}^{N-1}\overline{\varepsilon_0}(j)G_k^{\overline{(0,j)}},
\]
where the bar over $(0,j)$ means reduction modulo~$N$ (while $\overline{\varepsilon_0}$ means  the complex conjugate of~$\varepsilon_0$ as in Section \ref{ss:Gauss_sum_Bernoulli_Numbers}) and
\[
G_k^{\overline{(0,j)}}(z)=\sum\limits_{\substack{(c,d)\in\z^2\setminus\{(0,0)\}\\ (c,d)\equiv (0,j)\pmod{N}}}\frac{1}{(cz+d)^k}.
\]
Therefore for any $0\le j\le N-1$ coprime to $N$, we have
\[
\Big(\alpha_MG_k^{\overline{(0,j)}}\Big)\vert_k\gamma(z)=
\sum\limits_{\substack{(c,d)\in\z^2\setminus\{(0,0)\}\\ (c,d)\equiv (0,j)\pmod{N}}}\frac{1}{((cMu+dv)z+cM\beta+d\delta)^k}.
\]
Hence its constant term is given by
\[
\Upsilon_j=\sum\limits_{\substack{(c,d)\in\z^2\setminus\{(0,0)\}\\ (c,d)\equiv (0,j)\pmod{N} \\ cMu+dv=0}}
\frac{1}{(cM\beta+d\delta)^k}.
\]
Let $r=\gcd(v,M)$. Then $cMu+dv=0$ if and only if $cM'u+dv'=0$ where $M'=M/r$ and $v'=v/r$. 

If $u=0$, then $\Upsilon_j=0$ (for any $j$) unless $N=1$ in which case $j=0$ and 
\[
\Upsilon_0=\sum_{c\in\z}\frac{1}{(cM\beta)^k}=2\frac{1}{M^k}\zeta(k)
\]
since in this case $\beta=\pm1$ and $k$ is even. Therefore, according to Proposition~\ref{prop:FE_for_L_functions}, the constant term $\Upsilon$ of~$(\alpha_MG)\vert_k\gamma$ is $0$ if $N>1$ (and thus $N\nmid v=\pm1$) and is $\displaystyle{-\frac{C_k}{M^k}\frac{B_k}{k}}$ when $N=1$. Hence the result in this case.

Assume now $u\not=0$. Given $d\in\z$, $d\equiv j\pmod{N}$, the following conditions are then equivalent~:
\begin{enumerate}
\item there exists $c\in\z$, $c\equiv 0\pmod{N}$ such that $cMu+dv=0$;
\item we have $N\mid v'$ and $M'u\mid d$.
\end{enumerate}
Indeed, if the first condition is satisfied, we have $cM'u+dv'=0$ and thus $M'u\mid dv'$. But $M'u$ and $v'$ are coprime hence $M'u\mid d$. Besides, $0=cM'u+dv'\equiv dv'\pmod{N}$ and since $d\equiv j\pmod{N}$ and $\gcd(j,N)=1$, we get that $N\mid v'$.

On the other hand, if the second condition holds, put $c=-\displaystyle{\frac{d}{M'u}v'=-\frac{d}{Mu}v}$. Then $c\in\z$ satisfies $cMu+dv=0$ and since $N\mid v'$ by assumption, we get $c\equiv 0\pmod{N}$.

Moreover if these equivalent conditions are satisfied, then we have
\[
cM\beta+d\delta=\frac{1}{u}(cu\beta M+du\delta)=\frac{1}{u}(du\delta-dv\beta)=\frac{d}{u}.
\]
Therefore the constant term $\Upsilon$ of $(\alpha_MG)\vert_k\gamma$ is~$0$ when $N\nmid v'$ and is otherwise given by
\begin{align*}
\Upsilon
&=\sum\limits_{\substack{j=0\\ \gcd(j,N)=1}}^{N-1}\overline{\varepsilon_0}(j)\Upsilon_j
 =\sum\limits_{\substack{j=0\\ \gcd(j,N)=1}}^{N-1}\overline{\varepsilon_0}(j)
\sum\limits_{\substack{d\in\z\setminus\{0\} \\ d\equiv j\pmod{N} \\ M'u\mid d}} 
\left(\frac{u}{d}\right)^k \\
& =\sum\limits_{\substack{j=0\\ \gcd(j,N)=1}}^{N-1}\overline{\varepsilon_0}(j)
\sum\limits_{\substack{d\in\z\setminus\{0\} \\ d\equiv j/(M'u)\pmod{N}}} 
\frac{1}{(M'd)^k} \\
& = \frac{\overline{\varepsilon_0}(M'u)}{M'^k} \sum\limits_{\substack{j=0\\ \gcd(j,N)=1}}^{N-1}
\sum\limits_{\substack{d\in\z\setminus\{0\} \\ d\equiv j/(M'u)\pmod{N}}} 
\frac{\overline{\varepsilon_0}(d)}{d^k} \\
& = \frac{\overline{\varepsilon_0}(M'u)}{M'^k} 
\sum\limits_{d\in\z\setminus\{0\}} \frac{\overline{\varepsilon_0}(d)}{d^k} \\
& =2 \frac{\overline{\varepsilon_0}(M'u)}{M'^k} L(k,\overline{\varepsilon_0})
\end{align*}
as $\overline{\varepsilon_0}(-1)=(-1)^k$. Besides, since $N\mid v$ and $u\delta-v\beta=1$, we have $\overline{\varepsilon_0}(u)=\varepsilon_0(\delta)$. Using  Proposition \ref{prop:FE_for_L_functions}, we therefore find that in this case the constant term $\Upsilon$ of~$(\alpha_MG)\vert_k\gamma$ is non-zero and given by 
\[
\Upsilon=- \frac{2C_kW(\overline{\varepsilon_0})}{N^k}\frac{\overline{\varepsilon_0}(M')\varepsilon_0(\delta)}{M'^k}\frac{B_{k,\varepsilon_0}}{2k}.
\]
Hence the result in this case as well.
\endproof

\subsubsection{A useful Eisenstein series}\label{sss:ES_again}
In this paragraph we state some results about Eisenstein series that will be used in Section \ref{s:modularitecyclotomique}. Let $N$ be a squarefree integer and $k$ be an even integer~$\ge4$. We denote by~\(\mathcal{E}_k(\Gamma_0(N))\) the space spanned by the Eisenstein series of weight~\(k\) and level~\(\Gamma_0(N)\).

For each prime divisor~$p$ of~$N$, let $\delta_p\in\{1,p^{k-1}\}$. Put
\[
E(z)=\left[\prod_{\substack{p\mid N \\ p \textrm{ prime }}}(U_p-\delta_p\Id)\right]\alpha_NE_k(z)\in\mathcal{E}_k(\Gamma_0(N)),
\]
where~\(U_p\) (for \(p\) a prime divisor of~\(N\)) is the \(p\)-th Hecke operator acting on~\(\mathcal{E}_k(\Gamma_0(N))\). We remark that, even if it is not included in the notation, the function~$E$ does depend on the choice of the parameters~$\{\delta_p : p |N\}$. For a prime number~\(p\) not dividing~\(N\), we denote by~\(T_p\) the classical \(p\)-th Hecke operator acting on~\(\mathcal{E}_k(\Gamma_0(N))\) The following proposition summarizes the main properties of~$E$.
\begin{prop}\label{prop:Eisenstein}
The Eisenstein series $E$ is a normalized Hecke eigenform of level~$\Gamma_0(N)$ such that
\[
\left\{
\begin{array}{ll}
T_pE=(1+p^{k-1})E & \textrm{if }p\nmid N \\
U_{p}E=(p^{k-1}/\delta_p)E & \textrm{otherwise}.
\end{array}\right.
\]
We have that
\begin{equation}\label{eq:expanded_form_for_E}
E=\sum_{M\mid N}(-1)^{|M|}\delta_M\alpha_M E_k,
\end{equation}
where $|M|$ is the number of prime divisors of~$M$ and 
\[
\delta_M=\prod_{\substack{p\mid M \\ p\text{ prime}}}\delta_p.
\]
Let $s$ be a cusp of~$X_0(N)$. It is then $\Gamma_0(N)$-equivalent to~$1/v$ for some $v\mid N$ and the constant term of the Fourier expansion of~$E$ at the cusp~$s$ is 
\[
-\frac{B_k}{2k}\prod_{p\mid N}\left(1-\delta_p\left(\frac{\gcd(p,v)}{p}\right)^k\right),
\]
where the product runs over the prime divisors of~$N$.
\end{prop}
\proof The Eisenstein series $E_k$ is a normalized Hecke eigenform of level~$1$ with eigenvalue $1+p^{k-1}$ for each prime number~$p$. Let $p$ be a prime not dividing~$N$. It then follows from the action of Hecke operators on~\(q\)-expansions of modular forms (see, for instance, \cite[Prop.~5.2.2]{DiSh05}) that \(T_p\alpha_NE_k=\alpha_NT_pE_k\). Since the Hecke algebra spanned by $\{ T_p : p \nmid N, U_p: p|N\}$  is commutative, we get $T_pE=(1+p^{k-1})E$. 

Let now $p$ be a prime dividing~$N$ and~\(M\) an integer dividing~\(N\). We have that
\begin{equation}\label{eq:Hecke_action_on_q_exp}
U_p\alpha_ME_k=\left\{
\begin{array}{ll}
\alpha_{M/p}E_k &\text{ if \(p\mid M\)} \\
(1+p^{k-1})\alpha_ME_k-p^{k-1}\alpha_{Mp}E_k & \text{otherwise}.
\end{array}
\right.
\end{equation}

Indeed, the case $p | M$ follows from \emph{loc. cit.}  (note however that Diamond and Shurman use the notation~\(T_p\) in the case~\(p\mid N\) as well) and the case $p \nmid M$ follows from a direct calculation using the~\(q\)-expansion~(\ref{eq:classical_Ek}) of~\(E_k\)..

Since~\(N\) is squarefree, we therefore get
\begin{align*}
U_p(U_p-\delta_p\Id)\alpha_NE_k  & =  (U_p^2-\delta_pU_p)\alpha_N E_k \\
 & = U_p\alpha_{N/p}E_k-\delta_p\alpha_{N/p}E_k \\
 & = (1+p^{k-1})\alpha_{N/p}E_k-p^{k-1}\alpha_N E_k-\delta_p\alpha_{N/p}E_k\\
 & = \frac{p^{k-1}}{\delta_p}\alpha_{N/p}E_k-p^{k-1}\alpha_N E_k \\
 & = \frac{p^{k-1}}{\delta_p}\left(\alpha_{N/p}E_k-\delta_p\alpha_N E_k\right)\\
 & = \frac{p^{k-1}}{\delta_p}(U_p-\delta_p\Id)\alpha_NE_k.
\end{align*}
Hence the result in this case as well.

Using Equation~(\ref{eq:Hecke_action_on_q_exp}), we easily prove by induction on the number of prime divisors of~\(N\) the `expanded form'~(\ref{eq:expanded_form_for_E}) for~\(E\).

Let $s$ be a cusp of~$X_0(N)$. Since $N$ is squarefree, we have that $s$  is $\Gamma_0(N)$-equivalent to~$1/v$ where $v\mid N$. Let $M$ be a divisor of~$N$. Then, according to Proposition ~\ref{prop:constant_term} the constant term of the Fourier expansion of~$\alpha_ME_k$ at the cusp~$s$ is 
\[
-\frac{B_k}{2k}\left(\frac{\gcd(v,M)}{M}\right)^k.
\]
Since~\(N\) is squarefree, the constant term of the Fourier expansion of~$E$ at the cusp~$1/v$ is then
\[
-\frac{B_k}{2k}\sum_{M\mid N}(-1)^{|M|}\delta_M\left(\frac{\gcd(v,M)}{M}\right)^k=-\frac{B_k}{2k}\prod_{\substack{p\mid N \\ p\textrm{ prime}}}\left(1-\delta_p\left(\frac{\gcd(p,v)}{p}\right)^k\right).
\]
This proves the result.
\endproof

%% file: modularitereductible.tex
\section{Modularity of odd reducible semisimple representations}\label{s:modularite_reductible}

The following theorem is presumably well-known to experts, but we provide  a proof due to lack of suitable reference.
\begin{thm}\label{thm:modularity_general}
Every odd representation which is the direct sum of two characters  arises from a cuspidal eigenform.  
\end{thm}
Let $\rho:\gal(\qbar/\q)\rightarrow\GL_2(\Flbar)$ be such a representation and fix a place $w$ of~$\qbar$ above $l$, as in the Introduction. Assume first that $\rho\simeq\trivial\oplus\varepsilon\chi_l^b$ where $\varepsilon$ is unramified at~$l$ and $0\le b\le l-2$. The oddness condition here means $\varepsilon(-1)=(-1)^{b+1}$. We now define two integers~$N\ge 1$ and $k\ge 2$ and a character $\varepsilon_0$.

Let us denote by~$N$ the Artin conductor of~$\varepsilon$. It is coprime to~$l$ by assumption. Moreover if we denote by $\varepsilon_0$ the Teichm\"uller lift (with respect to~$w$) of~$\varepsilon$ we may identify it with a primitive Dirichlet character of conductor~$N$. It satisfies $\varepsilon_0(-1)=(-1)^{b+1}$ unless $l=2$, in which case we have $\varepsilon_0(-1)=1$.

We define a `weight'~$k$ attached to~$\trivial\oplus\varepsilon\chi_l^b$ as follows~:
\[
k=\left\{
\begin{array}{ll}
4 & \textrm{if $b=0$ and $l=2$} \\
l & \textrm{if $b=0$ and $l\ge3$} \\
l+1 & \textrm{if $b=1$} \\
b+1 & \textrm{if $b\ge2$} \\
\end{array}
\right..
\]
Note that $\varepsilon_0(-1)=(-1)^k$ and $k-1\equiv b\pmod{l-1}$. Hence $\rho\simeq\trivial\oplus\varepsilon\chi_l^{k-1}$.
\begin{rem}\label{poids}
Serre's recipe (see~\cite[Eq.~(2.3.2)]{Ser87}) for the weight of such a representation is (with Edixhoven's notation, \cite{Edi92})
\[
k_{\rho}=\left\{
\begin{array}{ll}
l & \textrm{if $b=0$} \\
b+1 & \textrm{if $b\ge1$} \\
\end{array}
\right.
\]
while Edixhoven's definition gives $k(\rho)=b+1$ ({\em loc. cit.}). Our definition is motivated by the fact that we want to avoid working with Eisenstein series of weight~$1$ or~$2$ in the proof of Theorem \ref{thm:modularity_particular_case} below.
\end{rem}
Let $\lambda$ be the prime ideal induced by our fixed place~$w$ in the ring of integers of the number field generated by the values of~$\varepsilon_0$. We can now state a  special case of Theorem \ref{thm:modularity_general}.
\begin{thm}\label{thm:modularity_particular_case}
Let $(N,k,\varepsilon_0)$ as above. Then, the representation $\trivial\oplus\varepsilon\chi_l^b$ arises from an eigenform in~$\CuspF{k}{Np}{\varepsilon_0}$ for every prime number $p\nmid Nl$ such that $\lambda$ divides the non-zero algebraic number $\displaystyle{\frac{B_{k,\varepsilon_0}}{2k}(\varepsilon_0(p)p^k-1)}$.
\end{thm}
\proof Let $p$ be a prime number not dividing~$Nl$. We consider
\[
E=E_k^{\trivial,\varepsilon_0}-\alpha_pE_k^{\trivial,\varepsilon_0}
\]
where $E_k^{\trivial,\varepsilon_0}$ is the Eisenstein series defined in Equation~(\ref{eq:Eisenstein_series_main}) and compute the constant term, say $a_{\gamma}$, of the $q$-expansion of $E\vert_k\gamma$ for any $\gamma=
\begin{pmatrix}
u & \beta \\
v & \delta \\
\end{pmatrix}\in\SL_2(\z)$. By construction we have $\varepsilon_0(-1)=(-1)^k$. According to Proposition~\ref{prop:constant_term}, if \(N\nmid v\), then \(a_{\gamma}=0\). Besides, if \(Np\mid v\), then~\(\gamma\in\Gamma_0(Np)\) and \(E\vert_k\gamma=\varepsilon_0(\delta)E\). Hence \(a_{\gamma}=0\) in this case as well. Finally, if $N\mid v$ and $p\nmid v$, then Proposition~\ref{prop:constant_term} gives
\[
a_{\gamma}=-\varepsilon_0(\delta)\frac{B_{k,\varepsilon_0}}{2k}\left(1-\frac{\overline{\varepsilon_0}(p)}{p^k}\right).
\]

Therefore, under the assumption that $\lambda$ divides $(\varepsilon_0(p)p^k-1)B_{k,\varepsilon_0}/2k$, the reduction of $E$ modulo~$\lambda$ is a cuspidal eigenform with coefficients in~$\Flbar$ and eigenvalues $1+\varepsilon(q)q^{k-1}$ for every prime~$q\nmid Np$. According to~\cite[Lem.~6.11]{DeSe74} we can find a form $f\in\CuspF{k}{Np}{\varepsilon_0}$ which is an eigenform for the Hecke operators~$\{T_q\}_{q\nmid Np}$ with corresponding eigenvalues~$\{a_q\}_{q\nmid Np}$ satisfying (for $q\not=l$)
\(\widetilde{a_q}= 1+\varepsilon(q)q^{k-1}\). Hence the representation~$\trivial\oplus\varepsilon\chi_l^b=\trivial\oplus\varepsilon\chi_l^{k-1}$ arises from an eigenform in~$\CuspF{k}{Np}{\varepsilon_0}$ as claimed.
\endproof

\noindent {\bf Proof of Theorem }~\ref{thm:modularity_general}: let $\rho$ be an odd representation which is the direct sum of two characters, say $\rho=\nu\oplus\nu'$. We thus have $(\nu\nu')(-1)=-1$. Put $\mu=\nu^{-1}\nu'$ and write $\mu=\varepsilon\chi_l^b$ where $\varepsilon$ is unramified at~$l$ and $0\le b\le l-2$. It satisfies $\mu(-1)=-1$ (or, equivalently $\varepsilon(-1)=(-1)^{b+1}$). Let $N$, $k$ and~$\varepsilon_0$ be the invariants as defined above attached to the representation $\trivial\oplus\mu=\trivial\oplus\varepsilon\chi_l^{b}$. By the previous theorem the representation $\trivial\oplus\mu$ arises from a Hecke eigenform~$f$ in~$\CuspF{k}{Np}{\varepsilon_0}$ for infinitely many primes~$p\nmid Nl$. Let us now consider some characteristic~$0$ lift~$\nu_0$ of~$\nu$ with finite image. It may be identified with some primitive Dirichlet character of modulus~$r$, say. Put $g=f\otimes\nu_0$. By a result of Shimura~\cite[Prop.~3.64]{Shi94}, we have that $g$ belongs to 
$\CuspF{k}{\lcm(Np,r^2,rN)}{\varepsilon_0\nu_0^2}$. Moreover $g$ is an eigenform for the Hecke operators outside $Npr$ (see~\cite[p.~34]{Rib77}) and its attached mod~$l$ Galois representation is isomorphic to~$\rho_f\otimes\nu\simeq\rho$. We thus have proved Theorem~\ref{thm:modularity_general}.

%% file: modularitecyclotomique.tex
\section{Non-optimal levels of $\rep$}\label{s:modularitecyclotomique}

The following proposition justifies our study of reducible representations of the special form $\trivial\oplus\chi_l^{k-1}$.

\begin{prop}\label{ss:reducible_mod_reps}
 Let $f$ be a newform of weight~$k\ge 2$, squarefree level~$N$ and trivial Nebentypus. For a prime~$l$, if  $l\ge k-1$ and $\rho_{f}$ is reducible, we have that $\rho_{f}$ is isomorphic to~$\rep$.
\end{prop}

\proof By assumption, $\rho_{f}$ is the direct sum of two characters $\nu_1$ and $\nu_2$. We decompose each of them  as a product of a character~$\varepsilon_i$ ($i=1,2$), which is unramified at~$l$, and some power of the cyclotomic character. The local description of the representation~$\rho_{f}$ at~$l$ (\cite[Thm.~2.5-2.6]{Edi92}) shows that, under the assumption~$l\ge k-1$, the exponents of the cyclotomic characters are~$0$ and~$k-1$. On the other hand, the product~\(\varepsilon_1\varepsilon_2\) is trivial (as the form~\(f\) has trivial Nebentypus) and by the squarefreeness assumption the characters~$\varepsilon_i$ are also unramified at the primes dividing~$N$ and hence trivial. \endproof

In the rest of  this section we fix an even integer~\(k\ge2\)  and focus our attention on the level-raising problem for (odd) representations of the form $\trivial\oplus\chi_l^{k-1}$ with $2\le k\le l-3$. With regards to the problem of classifying the squarefree non-optimal levels attached to $\trivial\oplus\chi_l^{k-1}$, we propose the following conjecture.

\begin{conj}\label{conj}
Let $k$ be an even integer~$\ge4$ and assume $l>k+1$. Then the representation~$\trivial\oplus\chi_l^{k-1}$ arises from a weight-$k$ newform of squarefree level~$N$, with  $l\nmid N$, and trivial Nebentypus if and only if at least one of the following conditions holds~:
\begin{enumerate}
\item we have $(p^k-1)(p^{k-2}-1)\equiv 0\pmod{l}$ for every prime number~$p$ dividing~$N$ and there exists a prime divisor~$p_0$ of~$N$ such that $p_0^{k}\equiv 1\pmod{l}$;
\item we have $p^{k-2}\equiv 1\pmod{l}$ for every prime number $p$ dividing~$N$ and $l$ divides the numerator of $B_k/k$, where $B_k$ denotes the $k$-th Bernoulli number.
\end{enumerate}
\end{conj}

We are able to show the direct implication in this conjecture (see Theorem~\ref{thm:easy_direction}). Concerning the reverse implication, we prove  a weaker statement (see Theorem~\ref{eigen}).  On the other hand, Theorem \ref{main} settles the conjecture in  the case where $N$ is  prime. The  case $N=1$ follows in one direction from a result of Ribet (\cite[Lem.~5.2]{Rib75}) and from Deligne-Serre's lifting lemma (\cite[Lem.~6.11]{DeSe74}) in the other (cf. Corollary \ref{cor:level_one_case} and Remark \ref{rem:level_one}).

We can prove that Conjecture \ref{conj} is actually equivalent to saying that the classical (necessary) level-raising condition at a prime (away from the level) is sufficient (cf. Section~\ref{ss:level_raising_condition}). In Section~\ref{ss:reverse_implication} we combine Theorems~\ref{thm:easy_direction} and~\ref{eigen} with a result of Diamond (\cite{Di91}) to prove the prime-level case of the conjecture. Using magma (\cite{Mag97}), we have computationally checked the validity of the conjecture for fixed weights and levels in various ranges. In particular it holds true for $k=4$ and $N<5000$ and for $6\le k<32$ and $N<50$.

\subsection{Necessary conditions}
In this paragraph we prove the following statement which corresponds to the direct implication in Conjecture~\ref{conj}.
\begin{thm}\label{thm:easy_direction}
Let $k$ be an even integer and $N$ be a non-negative squarefree integer. Assume $k\geq 4$, $l>k+1$ and $l\nmid N$. If the representation $\rep$ arises from a weight-$k$ newform of level~$N$ and trivial Nebentypus, then at least one of the following assertions holds~:
\begin{enumerate}
\item we have $(p^k-1)(p^{k-2}-1)\equiv 0\pmod{l}$ for every prime number~$p$ dividing~$N$ and there exists a prime divisor~$p_0$ of~$N$ such that $p_0^{k}\equiv 1\pmod{l}$;
\item we have $p^{k-2}\equiv 1\pmod{l}$ for every prime number $p$ dividing~$N$ and $l$ divides the numerator of $B_k/k$.
\end{enumerate}
\end{thm}
The proof splits into two steps. We first deduce some weaker conditions from the local description of modular representations at primes dividing exactly once the level. In a slightly different form, this was already done in a joint paper by Dieulefait and the first author (\cite{BiDi14}), but we briefly repeat the argument here for the sake of conciseness. We then strengthen these conditions using some (new) computations about Eisenstein series from Section \ref{sss:ES_again} to obtain Theorem~\ref{thm:easy_direction}.

Let $k,l$ and $N$ as in the theorem. Recall that we have fixed a place~\(w\) of~\(\qbar\) above~\(l\). Assume that $\rep$ arises from some newform $f$ of weight $k$ and level~$\Gamma_0(N)$. Let $p$ be a prime dividing $N$. By~\cite[Prop.~2.8]{LoWe12}, the restriction of $\rho_f$ to a decomposition group $D_p$ at $p$ is $\mu\chi_l^{k/2}\oplus\mu\chi_l^{k/2-1}$ where $\mu$ is the (at most) quadratic unramified character thats maps a Frobenius at~$p$ to the reduction modulo~\(w\) of~$a_p(f)/p^{k/2-1}$. Therefore, we have the following equality between sets of characters of $D_p$~: 
\[
\{\trivial,\chi_l^{k-1}\}=\{\mu\chi_l^{k/2},\mu\chi_l^{k/2-1}\}.
\]
\begin{enumerate}
\item Assume that $\trivial=\mu\chi_l^{k/2}$. Then, in particular $p^k\equiv 1\pmod{l}$.
\item Assume that $\trivial=\mu\chi_l^{k/2-1}$. Then, $\widetilde{a_p(f)}= 1$ and in particular, $p^{k-2}\equiv1\pmod{l}$. 
\end{enumerate}

Let us now assume that the first assertion of the theorem is not satisfied. According to the above discussion we have
\[
p^{k-2}\equiv 1\pmod{l};\quad p^k\not\equiv 1\pmod{l}\quad\textrm{and}\quad \widetilde{a_p(f)}=1
\]
for every prime $p$ dividing~$N$ and we must show that $l$ divides the numerator of~$B_k/k$. This will be achieved using a careful study of the constant term of the Fourier expansion at various cusps of a specific Eisenstein series which we now introduce.

Let us consider the Eisenstein series~$E$ as in Section \ref{sss:ES_again} with parameters $\delta_p=p^{k-1}$ for every prime $p\mid N$ and write \(E(z)=\sum_{n\geq0}a_n(E)q^n\), where \(q=e^{2i\pi z}\). Then by assumption and Proposition~\ref{prop:Eisenstein}, we have~:
\[
\widetilde{a_p(f)}= \widetilde{a_p(E)},\quad\textrm{for all primes }p\not=l.
\]
Besides, both $E$ and $f$ are normalized Hecke eigenforms. Therefore, we get
\[
\widetilde{a_n(f)}= \widetilde{a_n(E)},\quad\textrm{for all prime-to-$l$ integers $n$}.
\]
We denote by  $\widetilde{f}$ and $\widetilde{E}$ the reductions modulo~$w$ of~$f$ and~$E$ respectively.  Applying the  operator $\Theta=q \frac{d}{dq}$ (see \cite{Ser73}) we obtain the equality $\Theta(\widetilde{f})=\Theta(\widetilde{E})$. Since the $\Theta$ operator is injective under the assumption~$l>k+1$ (\cite[Cor.~3]{Kat77}), we conclude that $\widetilde{E} = \widetilde{f}$. In particular, $w$ divides the numerator of the constant term of $E$ at the cusp $\infty$. By Proposition~\ref{prop:Eisenstein}, this means that $l$ divides the numerator of $\displaystyle{\frac{B_k}{2k}\prod_{p\mid N}(1-p^{k-1})}$. Since $p^{k-1}\equiv 1\pmod{l}$ would imply $p^k\equiv1\pmod{l}$ (as $p^{k-2}\equiv1\pmod{l}$), contrary to the hypotheses, we get the desired result. This ends the proof of Theorem ~\ref{thm:easy_direction}.

\begin{rem}
According to Proposition~\ref{prop:Eisenstein}, the vanishing modulo~$l$ of the constant terms of~$E$ at the other cusps of~$X_0(N)$ does not give additional information.
\end{rem}

\subsection{Weaker converse statement and the prime level case}\label{ss:reverse_implication}
In what follows we present a weaker statement in the direction of the reverse implication of Conjecture~\ref{conj}.  We will finish this paragraph with a proof of Theorem~\ref{main}.
\begin{thm}\label{eigen}
Let $N$ be a positive squarefree integer. Let $k$ be an even integer~$\ge4$ and assume $l>k+1$. Assume that at least one of the following conditions holds~:
\begin{enumerate}
\item we have $(p^k-1)(p^{k-2}-1)\equiv 0\pmod{l}$ for every prime number~$p$ dividing~$N$ and there exists a prime divisor~$p_0$ of~$N$ such that $p_0^{k}\equiv 1\pmod{l}$;
\item we have $p^{k-2}\equiv 1\pmod{l}$ for every prime number $p$ dividing~$N$ and $l$ divides the numerator of $B_k/k$.
\end{enumerate}
Then the representation~$\trivial\oplus\chi_l^{k-1}$ arises from a weight-$k$ eigenform of  level~$N$ and trivial Nebentypus.
\end{thm}

\begin{rem}
This statement is weaker than the reverse implication in Conjecture \ref{conj} because it is not guaranteed that the eigenform is a newform.
\end{rem}

\proof Assume that either condition of the theorem is satisfied and let us consider the Eisenstein series~$E$ of Section~\ref{sss:ES_again} with the following choice of parameters~: 
\[
\delta_p=\left\{\begin{array}{ll}
                1 & \textrm{if }p^k\equiv1\pmod{l} \\
		p^{k-1} & \textrm{otherwise}
                \end{array}
\right..
\]
Recall from Equation~(\ref{eq:expanded_form_for_E}) that we have
\begin{equation}\label{eq:E_expanded}
E=\sum_{M\mid N}(-1)^{|M|}\delta_M\alpha_M E_k,
\end{equation}
where $|M|$ is the number of prime divisors of~$M$ and 
\[
\delta_M=\prod_{\substack{p\mid M \\ p\textrm{ prime}}}\delta_p.
\]
According to the assumptions $l>k+1$, $l\nmid N$ and the Van Staudt-Clausen theorem, the series $E$ has $l$-integral rational Fourier coefficients at~$\infty$. Let us denote by~$F$ its reduction modulo~$l$. It is a well-defined modular form over~$\FF_l$.

We now prove that $F$ is actually cuspidal. By Proposition~\ref{prop:Eisenstein}, this is clear under assumption (2) (as in particular, $l$ divides the numerator of $B_k/k$). Else, if we assume assumption (1), then there exists a prime divisor $p_0$ of~$N$ such that $\delta_{p_0}=1$. Let $s$ be a cusp of $X_0(N)$. It is $\Gamma_0(N)$-equivalent to some cusp of the form $1/v$ with $v\mid N$ and $1\le v\le N$. Then 
\[
1-\delta_{p_0}\left(\frac{\gcd(p_0,v)}{p_0}\right)^k=\left\{
\begin{array}{ll}
0 & \textrm{if }\gcd(p_0,v)=p_0 \\
1-p_0^{-k}& \textrm{ otherwise}.
\end{array}
\right.
\]
is congruent to~$0$ modulo~$l$. Hence the result by Proposition~\ref{prop:Eisenstein}. 

As it is already the case for~$E$, the cuspidal form~$F$ is a Hecke eigenform at level~$N$. Therefore according to the Deligne-Serre lifting lemma, there exist a finite extension $K/\q_l$ with ring of integers~$\mathcal{O}$ and uniformizer~$\mathcal{L}$ and a normalized Hecke eigenform $f\in\mathcal{S}_k(\Gamma_0(N);\mathcal{O})$ with system of eigenvalues~$\{c_p\}_p$ where $p$ runs over the primes, such that
\[
c_p\equiv 1+p^{k-1}\pmod{\mathcal{L}}\textrm{ if $p\nmid N$ and }c_p\equiv p^{k-1}/\delta_p\pmod{\mathcal{L}}\textrm{ otherwise}.
\]
Moreover $f$ is a classical modular form (as its Fourier coefficients are roots of the characteristic polynomials of the Hecke operators). 
\endproof
By a direct combination of Theorems~\ref{thm:easy_direction} and~\ref{eigen} we get a new proof of  the result  of Ribet mentioned in the Introduction, which constitutes the level~$1$ case of Conjecture~\ref{conj}. 
\begin{cor}\label{cor:level_one_case}
Let $k$ be an even integer~$\ge4$ and $l$ be a prime~$>k+1$. Then the representation~$\rep$ arises from a weight-$k$ eigenform of level~$1$ if and only if $l$ divides the numerator of~$B_k/k$.
\end{cor}
\begin{rem}\label{rem:level_one}
The direct implication of this result is due to Ribet \cite[Lem.~5.2]{Rib75}. The reverse implication is mentioned in \cite[Prop.~1]{Gha02}.
\end{rem}

\noindent {\bf Proof of Theorem~\ref{main}}: The  direct implication is a particular case of Theorem~\ref{thm:easy_direction}. 

Now we prove the reverse implication. By Theorem~\ref{eigen}, we have that $\rep$ arises from  an eigenform $f_0 \in\CuspFtriv{k}{N}$. If $f_0$ is a newform, then we are done. Hence, in what follows we will assume that $f_0$ is an oldform. We denote by $f$ its  associated (normalized) level~$1$ eigenform. By a standard application of the Cebotarev density theorem, the mod~$l$ representations $\rho_f$ and $\rho_{f_0}$ are isomorphic. 

Let $K$ be the number field spanned by the Fourier coefficients of $f$. Since $\rho_{f_0}$ is isomorphic to $\trivial\oplus\chi_l^{k-1}$ and $l\nmid N$, we have that there is an integral prime ideal $\lambda \subset O_K$ above~$l$ such that  $a_N(f)\equiv 1+N^{k-1} \pmod{\lambda}$.   We claim that 

\begin{equation}\label{DiamondHip}
a_N(f)^2\equiv N^{k-2}(1+N)^2 \pmod{\lambda}.
\end{equation}
Indeed,   
\[
a_N(f)\equiv 1+N^{k-1} \equiv 
\left\{\begin{array}{ll}
1+N^{-1} \pmod{\lambda} & \textrm{ if } N^k \equiv 1 \pmod{l} \\ 
1+N \pmod{\lambda} & \textrm{ if } N^{k-2} \equiv 1 \pmod{l}.
\end{array}\right.
\]
Since
\[
N^{k-2}(1+N)^2 \equiv 
\left\{\begin{array}{ll}
(1+N^{-1})^2 \pmod{l} & \textrm{ if } N^k \equiv 1 \pmod{l} \\ 
(1+N)^2 \pmod{l} & \textrm{ if } N^{k-2} \equiv 1 \pmod{l}
\end{array}\right.
\]
this proves the claim.

Relation~\eqref{DiamondHip} allows us to use a theorem of Diamond  (\cite[Thm.~1]{Di91}\footnote{Note that  Diamond's $(N,l,p)$ in {\em loc. cit.} is   $(1, N,l)$ in our notation.}), to ensure that there exists a normalized newform $f_1 \in S_k\big( \G(N)\big)^{\mathrm{new}}$ with eigenvalues in a finite extension $K'/K$ and an ideal $\lambda' \subset O_{K'}$ above $\lambda$ such that $$ a_p(f) \equiv a_p(f_1) \pmod{\lambda'} \textrm{ for all primes } p \nmid Nl.$$  Then, $\rho_{f_1}$ is isomorphic to $\rep$, concluding the proof.

\subsection{Relationship with the level-raising condition}\label{ss:level_raising_condition}
Let $$\rho~:\gal(\qbar/\q)\longrightarrow\GL_2(\Flbar)$$ be an odd semisimple  Galois representation of conductor~$N(\rho)$ (coprime to~\(l\)). We shall say that~$\rho$ satisfies the level-raising condition at a prime number~$p\nmid N(\rho)l$ if
\[
p\left(\tr\rho(\Frob_p)\right)^2=(1+p)^2\det\rho(\Frob_p)\quad\textrm{in }\Flbar,
\]
where $\Frob_p$ denotes a Frobenius element at~$p$ in~$\gal(\qbar/\q)$. Such a condition is satisfied if the representation~$\rho$ arises from a newform in~$\CuspFsans{k}{Np}$ with $(N,p)=1$. In particular, it is a necessary condition to raise the level of a modular representation from~$N$ to~$Np$. In their paper~\cite{DiTa94}, Diamond and Taylor prove that this is also sufficient when $\rho$ is assumed to be irreducible.

In the special case of the representation $\trivial\oplus\chi_l^{k-1}$, with even~$k\ge 2$, the level-raising condition at a prime~$p\not=l$ is merely
\begin{equation}\label{arriba}
 (p^k-1)(p^{k-2}-1)\equiv0\pmod{l}.
\end{equation}

If $k=2$, the congruence \eqref{arriba}  is automatically fulfilled for every~$p$, even though there are primes $p$ such that   $\trivial\oplus\chi_l$ is not  modular of weight~\(2\) and level~$p$. However, we believe  that the case $k\geq 4$ is  different and by analogy with the irreducible case, we propose the following conjecture:

\begin{conj}\label{conj2}
Let $k\ge 4$ be an even integer and $l$ be a prime~$>k+1$. Assume that $\trivial\oplus\chi_l^{k-1}$ arises from a newform in~$\CuspFtriv{k}{N}$ with $N$ squarefree and coprime to~$l$ and that $p\nmid Nl$ is a prime number at which $\trivial\oplus\chi_l^{k-1}$ satisfies the level raising condition, namely $(p^k-1)(p^{k-2}-1)\equiv0\pmod{l}$. Then, the representation $\trivial\oplus\chi_l^{k-1}$ arises from a newform in~$\CuspFtriv{k}{Np}$.
\end{conj}
Using Theorems~\ref{thm:easy_direction} and~\ref{eigen} above we now prove the following result.
\begin{prop}
Conjectures~\ref{conj} and~\ref{conj2} are equivalent.
\end{prop}
\proof Assume Conjecture \ref{conj} and the hypothesis of Conjecture~\ref{conj2}. Then Theorem~\ref{thm:easy_direction} ensures that the squarefree integer $Np$ satisfies the hypothesis of Conjecture \ref{conj}, thus proving Conjecture \ref{conj2}. 

Assume conversely Conjecture~\ref{conj2} and let us show that Conjecture~\ref{conj} holds. The direct implication therein corresponds to Theorem~\ref{thm:easy_direction}. Let us now prove the reverse implication. Let $N$ be a prime-to-$l$ squarefree integer that satisfies at least one of the conditions in the statement of Conjecture \ref{conj}. According to Theorem~\ref{eigen}, $\trivial\oplus\chi_l^{k-1}$ arises from a newform in~$\CuspFtriv{k}{M}$ for some integer~$M\mid N$. If $M\not=N$, we want to show that we can now raise the level from~$M$ to~$N$. Let $p$ be a prime dividing~$N/M$. Then one clearly has $(p^k-1)(p^{k-2}-1)\equiv0\pmod{l}$ and Conjecture~\ref{conj2} shows that $\trivial\oplus\chi_l^{k-1}$ arises from a newform in~$\CuspFtriv{k}{pM}$. If $pM=N$, we are done. Otherwise we can repeat this process until we reach~$N$. Hence, we have proved Conjecture~\ref{conj}.
\endproof 

\begin{rem}
As for irreducible representations, we deduce from Conjecture \ref{conj} that if the representation $\trivial\oplus\chi_l^{k-1}$ for~\(k\ge4\)  arises from newforms in~$\CuspFtriv{k}{M}$ and in~$\CuspFtriv{k}{N}$ for squarefree integers~$M\mid N$, then it also arises from a newform in~$\CuspFtriv{k}{N'}$ for any intermediate level~$M\mid N'\mid N$. As noticed by Ribet, such a result is false for $k=2$. The representation $\trivial\oplus\chi_5$ for instance arises in levels~$11$ and ~$66$ but neither in level~$22$ nor in level~$33$.
\end{rem}

%% file: corps.tex
\section{Lower bound for the highest degree of the coefficient field of  newforms}\label{corps}

In this section we prove Theorem \ref{thm:degree}. For a nonzero integer $m$, let $P^+(m)$ be the largest prime factor of $m$.  Let  
\begin{equation}\label{defp}
 \mathcal{P} := \{ N \textrm{ prime such that  } P^+(N-1) > N^{1/4}\}.
 \end{equation}
  That is, for every $N \in \mathcal{P}$,   there exists a prime $l$ with 
\begin{equation} \label{conditions}
 N \equiv 1 \pmod{l} \textrm{ and } l>N^{1/4}.
\end{equation} 

Let $A \subset \ene$ be a set consisting only of prime numbers. For $x \in \erre$, let 
$$A(x)=|\{a \in A : a \leq x\}|, \quad \pi(x)=|\{p \leq x : p \textrm{ is prime}\}|.$$ We recall that the quantity $$\liminf_{x \rightarrow \infty} \frac{A(x)}{\pi(x)}$$ is called the natural lower density of $A$. 

\begin{lema}
The set $\mathcal{P}$ has natural lower density at least $3/4$.
\end{lema}
\proof  In \cite{Luca},  Theorem 1, it is proved that  
\begin{equation} \label{densite}
\mathcal{P}(x) \geq \frac{3}{4} \cdot \frac{x}{\log x} + O \Big ( \frac{x}{(\log x)^{5/3}} \Big), \quad \textrm{as }x \rightarrow \infty.
\end{equation}
Using \eqref{densite} and the prime number theorem we obtain
$$\liminf_{x \rightarrow \infty} \frac{S(x)}{ \pi(x) } \geq \frac{3}{4},$$
as desired.
\endproof 
\noindent {\bf Proof of Theorem \ref{thm:degree}}: Let $\mathcal{P}$ be defined by \eqref{defp}. Take $N \in \mathcal{P}$ and a prime~$l$ as in  \eqref{conditions}. Assume $N \geq (k+1)^4$.   Then $N^k \equiv 1 \pmod{l}$ and $l>k+1$. Hence,  Theorem~\ref{main} and Mazur's Theorem \cite[Prop.~(5.12)]{Maz77}) ensure  that $\rep$ arises from a newform~$f=q+\sum_{n\ge2}a_nq^n$ of trivial Nebentypus, level~$N$ and weight $k$ if $k \geq 4$ and $k=2$ respectively. Put  $K=K_f$ and $d=[K_f:\q]$. Take a prime ideal $\lambda \subset O_K$ with $\lambda | l$ such that $$a_p \equiv 1+p^{k-1}\pmod{\lambda}, \quad \textrm{ for all primes } p \nmid Nl.$$

Moreover, $Nl$ is odd, so that we may consider this congruence for $p=2$. Deligne's bound (\cite{De74}, Th\'eor\`eme 8.2) implies that, for every archimedean place $\tau$ of $K$, we have that $|\tau (a_2)|\leq 2 \cdot 2^{(k-1)/2}$. Hence, the algebraic integer $b:=a_2- 1-2^{k-1} \in O_K$ is nonzero, belongs to $\lambda$ and satisfies $|\tau(b)|\leq (1+2^{(k-1)/2})^2$ for all $\tau$ as before. In particular, we have that 
$$ |N_{K/\q} (b)| \leq \left(1+2^{(k-1)/2}\right)^{2d}.  $$ Since $l|N_{K/\q} (b)$, we conclude that $l \leq   (1+2^ {(k-1)/2})^{2d}$, implying

$$d^{\mathrm{new}}_k(N) \geq d \geq \frac{\log l}{2\log  (1+2^{(k-1)/2}) } \geq   \frac{\log N}{8\log  (1+2^{(k-1)/2}) }.$$
This ends the proof of Theorem~\ref{thm:degree}. 

 \subsection{Final remarks}
The basic idea of using $a_2$  comes from  the proof of a similar statement in weight~$2$ by Dieulefait, Jimenez Urroz and Ribet (\cite[\S2]{DiUrRi11}). We are able to obtain a more general result because of  our Theorem~\ref{main} (that generalizes Mazur's theorem to higher weight) and the information on primes $p$ with large prime factors of $p-1$ given by Theorem~1 from~\cite{Luca}.

It is conjectured that for any $\eps>0$, the set prime numbers $p$ such that $P^+(p-1) \geq p^{1-\eps}$ has a positive lower density $\kappa(\eps) >0$. The bound  $\kappa(3/4 ) \geq 3/4$  is established   in \cite{Luca}  by extending a method of  Goldfeld who had previously obtained  $\kappa(1/2 ) \geq 1/2$ (\cite{Gold69}). Much effort has been invested in solving this conjecture for values of $\eps$ as small as possible (cf. \cite{Fou85}, \cite{HaBa96}). For  the purposes of Theorem \ref{thm:degree}, progress in this difficult problem would improve on the value of the constant $c_k$. However, such improvements would not change the fact that our method produces a constant $c_k$ that tends to zero with $k$. 

On the other hand, any value of $\eps$ \emph{bigger than} $3/4$ for which one could prove $\kappa(\eps) > 3/4$ would enlarge the set of primes for which our bound is valid (at the expense of a small loss in the constant $c_k$), thus improving Theorem~\ref{thm:degree} in an interesting way. For a nice compilation of conjectures and results about the density of this and related sets, see Section~2 of~\cite{BFPS04}.

If we assume Conjecture~\ref{conj}, it is possible to show an analogous lower bound for $d^{\mathrm{new}}_k(N)$ for $N$ in an appropriate family of squarefree integers.  Let $r$ be a non-negative integer. Put
\[
\mathcal{N}_r = \left\{ N \in \ene : N =p_1p_2\cdots p_r,  \omega(N)=r,  P^+\left( \gcd\limits_{1\le i\le r}(p_i-1) \right) > N^{\frac{1}{2r}} \right\}, 
\]
where $\omega(m)$ is number of different prime factors of the integer $m$ and $p_1,\ldots,p_r$ denote primes. It is shown in \cite{Luca} Theorem 2, that, as $x \rightarrow \infty$, we have~:

$$\frac{x^{\frac{1}{2}+ \frac{1}{2r}}}{(\log x)^{r+1}} \ll_r |\{   N \in \mathcal{N}_r: N \leq x\}| \ll_r \frac{x^{\frac{1}{2}+ \frac{1}{2r}}(\log \log x)^{r-1}}{(\log x)^{2}}.$$

These estimates show that the set $\mathcal{N}_r$ is infinite and that,  if $r \geq 2$, this set has density zero when regarded as a subset of squarefree numbers with exactly $r$ prime divisors. 

Mimicking the argument given above when $N$ is prime, we finally prove the following result.

\begin{thm} \label{souslaconjecture}
Assume Conjecture~\ref{conj} and let~\(q\) be a fixed prime number. Then, for every integer $r\ge2$ and every even~\(k\geq 4\), we have that 
$$d^{\mathrm{new}}_k(N) \gg_{k,q} \frac{1}{r}\log N, \quad \text{as }N\rightarrow \infty, N \in \mathcal{N}_r\text{ coprime to~\(q\)}.$$
\end{thm}

\noindent \emph{Sketch of proof.} Consider~\(N\in\mathcal{N}_r\). By assumption, there exists a prime~\(l\) such that \(N\equiv1\pmod{l}\) and~\(l>N^{\frac{1}{2r}}\). Moreover, if~\(N\) is large enough, then~\(l>k+1\). By Conjecture~\ref{conj}, there exists a newform~\(f\in\CuspFtriv{k}{N}\) with eigenvalues~\(\{a_p\}_p\) giving rise to~\(\trivial\oplus\chi_l^{k-1}\), that is~:
\[
a_p \equiv 1+p^{k-1}\pmod{\lambda}, \quad \textrm{ for all primes } p \nmid Nl,
\]
where~\(\lambda\) is some prime ideal in~\(\qbar\) above~\(l\).  Besides, \(N\) is coprime to~\(q\) by assumption and thus for large enough~\(N\) we have~\(q\nmid Nl\). We then conclude as in the proof of Theorem~\ref{thm:degree} using~\(a_q\) instead of~\(a_2\).
